\documentclass[12pt,leqno]{amsart}
\usepackage{amssymb}
\usepackage{amsmath,amssymb,color}
\numberwithin{equation}{section}

\newcommand{\weight}{e^{2s\alpha}}

\newcommand{\ep}{\varepsilon}
\newcommand{\la}{\lambda}
\newcommand{\va}{\varphi}
\newcommand{\ppp}{\partial}

\newcommand{\www}{\widetilde}

\newcommand{\R}{\mathbb{R}}

\newcommand{\N}{\mathbb{N}}

\newcommand{\ooo}{\overline}
\newcommand{\OOO}{\Omega}

\newcommand{\sumij}{\sum_{i,j=1}^d}

\newcommand{\HHHH}{H^{2,1}}
\newcommand{\EEEE}{\Vert \eta\Vert_{C(\overline{\Omega})}}

\newcommand{\hhalf}{\frac{1}{2}}

\newcommand{\ddd}{\mathrm{div}\,}  

\allowdisplaybreaks

\title
[]
{
Lipschitz stability for determination of states and 
inverse source problem for the mean field game equations
}

\author{$^1$ Oleg Imanuvilov, $^2$ Hongyu Liu and $^3$ Masahiro Yamamoto}

\thanks{
$^1$ Department of Mathematics, Colorado State University, 101 Weber Building, 
Fort Collins CO 80523-1874, USA,\\
e-mail: {\tt oleg@math.colostate.edu}
\\
$^2$ Department of Mathematics, City University of Hong Kong, Kowloon, 
Hong Kong SAR, China\\
email: {\tt hongyliu@cityu.edu.hk}
\\
$^3$ Graduate School of Mathematical Sciences, The University
of Tokyo, Komaba, Meguro, Tokyo 153-8914, Japan \\
e-mail: {\tt myama@ms.u-tokyo.ac.jp}
}

\date{}
\begin{document}
\maketitle

\begin{abstract}
We consider solutions satisfying the zero Neumann boundary condition and
a linearized mean field game equation in $\OOO \times (0,T)$ whose 
principal coefficients depend on the time and spatial variables
with general Hamiltonian, 
where $\OOO$ is a bounded domain in $\R^d$ and $(0,T)$ is the time
interval.  We first prove the Lipschitz stability in $\OOO \times (\ep, T-\ep)$
with given $\ep>0$ 
for the determination of the solutions by Dirichlet data on arbitrarily 
chosen subboundary of $\ppp\OOO$.
Next we prove the Lipschitz stability for an inverse problem of determining 
spatially varying factors of source terms and a coefficient
by extra boundary data and spatial data at intermediate time.    
\end{abstract} 
\baselineskip 18pt

\section{Introduction}

Recently the mean field game has called great attention and we refer
for example, to
Achdou, Cardaliaguet, Delarue, Porretta and Santambrogio
\cite{ACDPS}, Cardaliaguet, Cirat and Porretta \cite{CCP},
Lasry and Lions \cite{LL} and the references therein.

Let $\OOO \subset \R^d$ be a smooth bounded domain and let 
$\nu = \nu(x) = (\nu_1(x), ..., \nu_d(x))$ be the outward unit normal 
vector to $\ppp\OOO$ at 
$x \in \ppp\OOO$.  Let  
$$
Q:= \OOO \times (0,T).
$$
Then one typical mean field game system can be described by 
$$
\left\{ \begin{array}{rl}
& \ppp_tu + a(x,t)\Delta u(x,t)) 
- \hhalf\kappa(x,t)\vert \nabla u(x,t)\vert^2 + p(x,t)v = F(x,t), \\
& \ppp_tv - \Delta (a(x,t)v(x,t)) - \ddd(\kappa(x,t)v(x,t)\nabla u(x,t))
= G(x,t), \quad (x,t) \in Q
\end{array}\right.
                                                          \eqno{(1.1)}
$$
with the homogeneous Neumann boundary condition
$$
\nabla u \cdot \nu = \nabla v\cdot \nu = 0 
\quad \mbox{on $\ppp\OOO \times (0,T)$},                          \eqno{(1.2)}
$$
for example, and given $u(\cdot,T)$ and $v(\cdot,0)$ in $\OOO$.

Although there are many works on the well-posedness and related properties 
of solutions $u, v$ for the above system,
very few efforts have been devoted to inverse problems for (1.1).
We can refer to  Klibanov \cite{Kl23}, Klibanov and Averboukh \cite{KlAv},
Klibanov, Li and Liu \cite{KLL1}, \cite{KLL2},
Liu and Zhang \cite{LZ1}, \cite{LZ2}.

For inverse problems, we can mention two types for example:
\\
{\bf (i) Determination of state:} determine $u, v$ by extra data on a 
subboundary 
of $\ppp\OOO$ over a time interval.
\\
{\bf (ii) Inverse source problem:} Determine spatially varying 
factors of the source terms $F(x,t)$ and $G(x,t)$ by lateral data and spatial
data of $u, v$.
\\

As for the type (i), we can refer to \cite{KLL2}, and as for other types of
state determination for $u, v$ with data chosen among
$u(\cdot,0), u(\cdot,T), v(\cdot,0), v(\cdot,T)$, see 
\cite{Kl23}, \cite{KlAv}, \cite{KLL1}.  Also see Liu and Yamamoto \cite{LY} 
as related state determination problems.

However, to the best knowledge of the authors, there are no publications on
the uniqueness and the stability for the inverse source problems.
The main purpose of this article is to establish the Lipschitz stability
for the above two types of inverse problems, which have not been found in 
the existing articles.

We mainly consider a linearized equation of (1.1), which is formulated 
as follows.

We set 
$\ppp_i = \frac{\ppp}{\ppp x_i}$, $1\le i \le d$ and 
$\ppp_t = \frac{\ppp}{\ppp t}$ and for $\gamma := (\gamma_1, ...,
\gamma_d) \in (\N \cup \{ 0\})^d$ we define 
$\ppp_x^{\gamma}: = \ppp_1^{\gamma_1}\cdots \ppp_d^{\gamma_d}$ and
$\vert \gamma\vert := \gamma_1 + \cdots + \gamma_d$.

We introduce partial differential operators of the second order
with $(x,t)$- dependent 
coefficients by
$$
\left\{ \begin{array}{rl}
&A(t)u:= \sumij a_{ij}(x,t)\ppp_i\ppp_j u + \sum_{i=1}^d a_j(x,t)\ppp_ju
+ a_0(x,t)u, \\
&B(t)v:= \sumij b_{ij}(x,t)\ppp_i\ppp_j v + \sum_{i=1}^d b_j(x,t)\ppp_jv
+ b_0(x,t)v,\\
& A_0(t)u := \sum_{\vert \gamma\vert \le 2} b_{\gamma}(x,t) \ppp_x^{\gamma}u,
\end{array}\right.
$$
where $a_{ij}, b_{ij} \in C^1(\ooo{Q})$ and $a_{ij} = a_{ji}$,
$b_{ij} = b_{ji}$ for $1\le i,j\le d$, 
$c_0, a_k, b_k \in L^{\infty}(Q)$
for $0\le k \le d$, $b_{\gamma}\in L^{\infty}(Q)$ for 
$\vert\gamma \vert \le 2$, and it is assumed that there exists a constant 
$\chi > 0$ such that 
$$
\sumij a_{ij}(x,t)\xi_i\xi_j \ge \chi \sum_{j=1}^d \xi_j^2, \quad
\sumij b_{ij}(x,t)\xi_i\xi_j \ge \chi \sum_{j=1}^d \xi_j^2
\quad \mbox{for all $(x,t) \in Q$ and $\xi_1, ..., \xi_d \in \R$}.
$$
We define
$$
\ppp_{\nu_A}u := \sumij a_{ij}(\ppp_ju)\nu_i, \quad
\ppp_{\nu_B}v := \sumij b_{ij}(\ppp_jv)\nu_i  \quad
\mbox{on $\ppp\OOO\times (0,T)$}.
$$

We consider a linearized mean field game system:
$$
\left\{ \begin{array}{rl}
& \ppp_tu + A(t)u = c_0(x,t)v + F(x,t), \\
& \ppp_tv - B(t)v = A_0(t)u + G(x,t), \quad (x,t) \in Q
\end{array}\right.
                                          \eqno{(1.3)}
$$
with Neumann boundary condition
$$
\ppp_{\nu_A}u = \rho, \quad \ppp_{\nu_B}v = \sigma 
\quad \mbox{on $\ppp\OOO \times 
(0,T)$}.                                    \eqno{(1.4)}
$$

We emphasize that general Hamiltonians can be considered, which produces  
the second-order partial differential operator 
$A_0(t)u$ after the linearization.

In this article, we establish the Lipschitz stability results
for the following inverse problems.  Let $\Gamma$ be an arbitrarily chosen 
non-empty subboundary of $\ppp\OOO$, $t_0 \in (0,T)$ be arbitrarily given, and
let $u, v$ satisfy (1.3) and (1.4).
\\
{\bf Determination of state:}
{\it Determine $u, v$ in $Q$ by $u\vert_{\Gamma\times (0,T)}$ and
$v\vert_{\Gamma\times (0,T)}$.
}
\\
{\bf Inverse source problem:}
{\it
In (1.1) let $F(x,t) = q_1(x,t)f(x)$ and $G(x,t) = q_2(x,t)g(x)$ for
$(x,t) \in Q$.  Then determine $f, g$ in $\OOO$ by data
$u\vert_{\Gamma\times (0,T)}$, $v\vert_{\Gamma\times (0,T)}$
and $\{ u(\cdot,t_0), \, v(\cdot,t_0)\}$ in $\OOO$.
}
\\

We set
$$
\Vert u\Vert_{H^{2,1}(Q)}:= \left( \Vert u\Vert^2_{L^2(Q)}
+ \Vert \nabla u\Vert^2_{L^2(Q)} + \sumij \Vert \ppp_i\ppp_j u\Vert^2_{L^2(Q)}
+ \Vert \ppp_tu\Vert^2_{L^2(Q)}\right)^{\frac{1}{2}}.
$$

We first describe our main result for the state determination.
\\
{\bf Theorem 1.}
\\
{\it
We assume that $u,v \in H^{2,1}(Q)$ satisfy (1.3) and (1.4).
For arbitrarily given $\ep > 0$, we can find a constant $C_{\ep} >0$ such that 
\begin{align*}
& \Vert u\Vert_{\HHHH(\OOO\times (\ep,T-\ep))}
+ \Vert v\Vert_{\HHHH(\OOO\times (\ep,T-\ep))}
\le C_{\ep}(\Vert F\Vert_{L^2(Q)} + \Vert G\Vert_{L^2(Q)}\\
+ & \Vert u\Vert_{H^1(0,T;L^2(\Gamma))}
+ \Vert \nabla u\Vert_{L^2(\Gamma \times (0,T))}
+ \Vert v\Vert_{H^1(0,T;L^2(\Gamma))}
+ \Vert \nabla v\Vert_{L^2(\Gamma \times (0,T))}\\
+ & \Vert \rho\Vert_{H^1(0,T;H^{\hhalf}(\ppp\OOO))}
+  \Vert \sigma\Vert_{H^1(0,T;H^{\hhalf}(\ppp\OOO))}).
\end{align*}
}

In particular, we directly see
\begin{align*}
& \Vert u(\cdot,t)\Vert_{L^2(\OOO)}
+ \Vert v(\cdot,t)\Vert_{L^2(\OOO)}
\le C_{\ep}(\Vert F\Vert_{L^2(Q)} + \Vert G\Vert_{L^2(Q)}\\
+ & \Vert u\Vert_{H^1(0,T;L^2(\Gamma))}
+ \Vert \nabla u\Vert_{L^2(\Gamma \times (0,T))}
+ \Vert v\Vert_{H^1(0,T;L^2(\Gamma))}
+ \Vert \nabla v\Vert_{L^2(\Gamma \times (0,T))}\\
+& \Vert \rho\Vert_{H^1(0,T;H^{\hhalf}(\ppp\OOO))}
+  \Vert \sigma\Vert_{H^1(0,T;H^{\hhalf}(\ppp\OOO))}), \quad
\ep \le t \le T-\ep.
\end{align*}

We emphasize that Theorem 1 asserts unconditional stability which implies
that we do not need to impose any boundedness conditions for $u$ and $v$.

The norm of boundary values $\rho$ and $\sigma$ is not the best possible, but 
improvements require more arguments and here we do not pursue.
We can obtain a similar estimate for the case where $u, v$ are given on 
$\ppp\OOO \times (0,T)$ and $\ppp_{\nu_A}u, \ppp_{\nu_B}v$ on 
$\Gamma \times (0,T)$.
 
In Klibanov, Li and Liu \cite{KLL2}, 
the H\"older stability is proved with data $u, v, \nabla u, \nabla v$ 
on the whole lateral boundary $\ppp\OOO \times 
(0,T)$, which is the case of $\Gamma = \ppp\OOO$ in Theorem 1.
By the parabolicity of the equations (1.3), extra Dirichlet data should be 
limited to any small subboundary.
\\

Next we state our main result on the inverse source problem.
In (1.3) we assume 
$$
F(x,t) = q_1(x,t)f(x), \quad G(x,t) = q_2(x,t)g(x), \quad (x,t) \in Q,
$$
where $q_1, q_2 \in W^{1,\infty}(0,T;L^{\infty}(\OOO))$.
Moreover we assume
\begin{align*}
& a_{ij}, \ppp_ta_{ij} \in C^1(\ooo{Q}), \quad
c_0, \ppp_tc_0, a_k, \ppp_ta_k, b_k, \ppp_tb_k \in L^{\infty}(Q)
\quad \mbox{for $1\le i,j\le d$ and $0\le k \le d$}, \\
&b_{\gamma}, \ppp_tb_{\gamma} \in L^{\infty}(Q) \quad 
\mbox{for $\vert \gamma \vert \le 2$}.
\end{align*}
We fix $t_0 \in (0,T)$ arbitrarily and we are given
$$
u_0(x) := u(x,t_0), \quad v_0(x):= v(x,t_0), \quad x\in \OOO.
$$
Then
\\
{\bf Theorem 2 (global unconditional Lipschitz stability for an inverse 
source problem).}
\\
{\it
Let $u, v \in H^{2,1}(Q)$ satisfy (1.3), $\ppp_tu, \ppp_tv 
\in H^{2,1}(Q)$ and
$$
\ppp_{\nu_A}u = \ppp_{\nu_B}v = 0 \quad \mbox{on $\ppp\OOO \times
(0,T)$}.
$$
We assume
$$
\vert q_1(x,t_0)\vert, \, \vert q_2(x,t_0)\vert > 0 \quad \mbox{for all 
$x \in \ooo{\OOO}$}.
                                                  \eqno{(1.5)}
$$
Then there exists a constant $C>0$ such that 
\begin{align*}
& \Vert f\Vert_{L^2(\OOO)} + \Vert g\Vert_{L^2(\OOO)}\\
\le& C\left(\Vert u(\cdot,t_0)\Vert_{H^2(\OOO)}
+ \Vert v(\cdot,t_0)\Vert_{H^2(\OOO)}
+ \sum_{k=0}^1 (D(\ppp_t^ku)^2 + D(\ppp_t^kv)^2) \right).
\end{align*}
}

Here and henceforth we set 
$$
D(u) := \left(\int_{\Gamma\times (0,T)} (\vert \ppp_tu\vert^2 
+ \vert \nabla u\vert^2 + \vert u\vert^2) dSdt\right)^{\hhalf}.
$$

We emphasize that in our stability, we do not assume neither data  
$u(\cdot,T)$ nor $v(\cdot,0)$ in $\OOO$, 
nor any a priori bounds on $u,v,f, g$, but we require data $u(\cdot,t_0)$ 
and $v(\cdot,t_0)$ over $\OOO$ at an intermediate time $t_0 \in (0,T)$.
Our stability can be understood unconditional in the sense that 
we do not need to assume any a priori boundedness conditions.
\\
\vspace{0.1cm}

Our key is a classical Carleman estimate for a single parabolic equation
with singular weight function by Imanuvilov \cite{Ima}.
The linearized mean field game equations (1.3) have two features:
\begin{itemize}
\item
The equation in $u$ is backward and the one in $v$ is forward.
\item
The equation in $v$ contains the second-order spatial derivatives 
$A_0(t)u$ of $u$.
\end{itemize}

The mixed forward and backward equations in (1.3),
makes the direct problem such as initial boundary 
value problem difficult, but thanks to the symmetry of the time
variable in the weight function, this does not matter for 
Carleman estimates, our main tool.

Usually systems with coupled principal parts cause
difficulty for establishing relevant Carleman estimates.
However, in our case, although the second equation in (1.3) is coupled with 
the second-order terms of $u$, the first one in (1.3) is coupled only with 
zeroth order term of $v$, which enables us executing a typical argument
of absorbing the second-order terms of $u$ by taking large parameters
of the Carleman estimate.  

This article is composed of four sections.  In Section 2, we prove
a Carleman estimate (Theorem 3) for the linearized mean field game
equations and complete the proof of Theorem 1.  
Section 3 is devoted to the proof of Theorem 2.  In Section 4, we 
discuss the state determination problem for the original nonlinear 
mean value field equations (1.1).

\section{Carleman estimate for the mean field game equations and the proof 
of Theorem 1}
 
For an arbitrarily given subboundary $\Gamma \subset \ppp\OOO$, 
it is known (Lemma 2.3 in Imanuvilov and Yamamoto \cite{IY98}, 
Imanuvilov \cite{Ima}) that there exists 
$\eta \in C^2(\ooo{\OOO})$ such that 
$$
\eta > 0 \quad \mbox{in $\OOO$}, \quad \vert \nabla\eta \vert > 0
\quad \mbox{on $\ooo{\OOO}$}, \quad 
\ppp_{\nu_A}\eta\le 0, \quad \ppp_{\nu_B}\eta \le 0 
\quad \mbox{on $\ppp\OOO\setminus \Gamma$}.
$$
We set
$$
\va(x,t) = \frac{e^{\la\eta(x)}}{t(T-t)}, \quad
\alpha(x,t) = \frac{e^{\la\eta(x)} - e^{EE}}{t(T-t)}
$$
and
$$
t_0 = \frac{T}{2}.
$$
We set
$$
 M:= \sumij (\Vert a_{ij}\Vert_{C^1(\ooo{Q})}
+ \Vert b_{ij}\Vert_{C^1(\ooo{Q})})
+ \sum_{k=0}^d (\Vert a_k\Vert_{L^{\infty}(Q)}
+ \Vert b_k\Vert_{L^{\infty}(Q)}) + \sum_{\vert \gamma\vert\le 2}
\Vert b_{\gamma}\Vert_{L^{\infty}(Q)}.
$$

We recall 
$$
D(u) := \left( \int_{\Gamma\times (0,T)} 
(\vert \ppp_tu\vert^2 + \vert \nabla u\vert^2 
+ \vert u\vert^2) dSdt\right)^{\hhalf}.
$$

Henceforth $C>0$ denotes generic constants which are independent of $s, \la>0$
which are parameters, and we write $C(s,\la)$ when we need to specify the 
dependency.

Then we can prove
\\
{\bf Lemma 1.}
\\
{\it
There exists a constant $\la_0>0$ such that for each $\la > \la_0$, we can find
a constant $s_0(\la) > 0$ satisfying: there exists a constant $C>0$ such that 
\begin{align*}
& \int_Q \left(\vert \ppp_tu\vert^2 + \sumij \vert \ppp_i\ppp_ju\vert^2
+ s^2\la^2\va^2\vert \nabla u\vert^2 + s^4\la^4\va^4\vert u\vert^2\right)
\weight dxdt\\
\le& C\int_Q s\va\vert \ppp_tu + A(t)u\vert^2 \weight dxdt
+ C_1(\la,s)D(u)^2
\end{align*}
for all $s > s_0(\la)$ and $u\in H^{2,1}(Q)$ satisfying
$\ppp_{\nu_A}u = 0$ on $\ppp\OOO\times (0,T)$.
}
\\

Here the constant $C>0$ depends continuously on $M$: bound of the
coefficients and $\la$ but
independent of $s \ge s_0(\la)$, while $\la_0>0$ depends continuously
on $M$.
\\ 
{\bf Proof.}
By Fursikov and Imanuvilov \cite{FI} or Imanuvilov \cite{Ima}, we
know
\begin{align*}
& \int_Q \biggl\{ \frac{1}{s\va}\left(
\vert \ppp_tw\vert^2 + \sumij \vert \ppp_i\ppp_jw\vert^2\right)
+ s\la^2\va\vert \nabla w\vert^2 + s^3\la^4\va^3\vert w\vert^2 \biggr\}
\weight dxdt\\
\le& C\int_Q s\va\vert \ppp_tw - A(t)w\vert^2 \weight dxdt
+ C_1(\la,s)D(w)^2
\end{align*}
for $w \in H^{2,1}(Q)$ satisfying $\ppp_{\nu_A}w = 0$ on 
$\ppp\OOO\times (0,T)$.

We set $u(x,t) = w(x,T-t)$ for $(x,t) \in Q$.
Then, since $\alpha(x,t) = \alpha(x,T-t)$ and
$\va(x,t) = \va(x,T-t)$ for $(x,t) \in Q$, the change of varables 
$t \mapsto T-t$ verifies 
\begin{align*}
& \int_Q \biggl\{ \frac{1}{s\va}\left(
\vert \ppp_tu\vert^2 + \sumij \vert \ppp_i\ppp_ju\vert^2\right)
+ s\la^2\va\vert \nabla u\vert^2 + s^3\la^4\va^3\vert u\vert^2 \biggr\}
\weight dxdt\\
\le& C\int_Q \vert \ppp_tu + A(t)u\vert^2 \weight dxdt
+ C_1(\la,s)D(u)^2.
\end{align*}

Next we have to prove the estimate with one more factor $s\va$.
For simplicity we consider $A(t) = \Delta$.

By $\eta(x) \ge 1$ for $x\in \ooo{\OOO}$, we have
$$
\ppp_i\va = \la(\ppp_i\eta)\va, \quad
\ppp_i^2\va = (\la\ppp_i^2\eta + \la^2(\ppp_i\eta)^2)\va,
$$
and
$$
\frac{d}{dt}\left( \frac{1}{t(T-t)}\right) 
= \frac{2(t-t_0)}{t^2(T-t)^2}.
$$
Hence 
$$
\vert \ppp_t\va\vert \le C\va^2, \quad 
\vert \nabla\va\vert^2 \le C\la\va, \quad
\vert \ppp_i^2\va\vert \le C\la^2\va \quad \mbox{in $Q$},
$$
where the constant $C>0$ is independent of $\la > 0$.

We set $w:= \va^{\hhalf}u$.  Then
\begin{align*}
& \ppp_tw = \hhalf\va^{-\hhalf}(\ppp_t\va)u + \va^{\hhalf}\ppp_tu,\\
& \ppp_iw = \hhalf\va^{-\hhalf}(\ppp_i\va)u + \va^{\hhalf}\ppp_iu,\\
& \ppp_i\ppp_jw = \va^{\hhalf}\ppp_i\ppp_ju + \hhalf \va^{-\hhalf}
((\ppp_i\va)\ppp_ju + (\ppp_j\va)\ppp_iu)
+ \left( -\frac{1}{4}\va^{-\frac{3}{2}}(\ppp_i\va)\ppp_j\va
+ \hhalf \va^{-\hhalf}(\ppp_i\ppp_j\va)\right)u,
\end{align*}
and so
$$
\vert \ppp_iw\vert \le C(\la\va^{\hhalf}\vert u\vert 
+ \va^{\hhalf}\vert \ppp_iu\vert)
$$
and
\begin{align*}
& \ppp_tw + \Delta w = \va^{\hhalf}(\ppp_tu+\Delta u)
+ \hhalf\va^{-\hhalf}(\ppp_t\va)u\\
+ &\left( \hhalf \la\Delta \eta + \frac{\la^2}{4}\vert \nabla\eta\vert^2
\right)\va^{\hhalf}u + \la\va^{\hhalf}\nabla\eta \cdot \nabla u \quad
\mbox{in $Q$}.
\end{align*}
Hence,
$$
\vert \ppp_tw + \Delta w\vert \le C(\va\vert F\vert^2 + \va^3\vert u\vert^2
+ \la^4\va\vert u\vert^2 + \va\vert \nabla u\vert^2).
$$
Consequently the already proved Carleman estimate (\cite{Ima}) yields
$$
 \int_Q \biggl\{ \frac{1}{\va}\left(
\vert \ppp_tw\vert^2 + \sumij \vert \ppp_i\ppp_jw\vert^2\right)
+ s^2\la^2\va\vert \nabla w\vert^2 + s^4\la^4\va^3\vert w\vert^2 \biggr\}
\weight dxdt
$$
$$
\le C\int_Q s\va\vert F\vert^2 \weight dxdt
+ C\int_Q (s\va^3+s\la^4\va) \vert u\vert^2 \weight dxdt
+ C\int_Q s\va \vert u\vert^2 \weight dxdt 
+ CD(w)^2.                      \eqno{(2.1)}
$$
On the other hand,
$$
\left\{ \begin{array}{rl}
& \va \vert \ppp_tu\vert^2 \le C(\vert \ppp_tw\vert^2
+ \va^3\vert u\vert^2), \\
& \va\vert \nabla u\vert^2 \le C(\vert \nabla w\vert^2 
+ \la^2\va\vert u\vert^2), \\
& \va\vert \ppp_i\ppp_ju\vert^2
\le C(\vert \ppp_i\ppp_jw\vert^2 + \la^2\va\vert \nabla u\vert^2
+ \la^4\va\vert u\vert^2),
\end{array}\right.
                             \eqno{(2.2)}
$$
so that
\begin{align*}
& \vert \ppp_tu\vert^2 + \sumij \vert \ppp_i\ppp_ju\vert^2\\
\le& \frac{C}{\va} \left( \vert \ppp_tw\vert^2 
+ \sumij \vert \ppp_i\ppp_jw\vert^2\right) + C(\la^2\vert \nabla u\vert^2
+ \la^4\vert u\vert^2 + \va^2 \vert u\vert^2).
\end{align*}
Therefore, by (2.2) we obtain
\begin{align*}
& \vert \ppp_tu\vert^2 + \sumij \vert \ppp_i\ppp_ju\vert^2
+ s^2\la^2\va^2 \vert \nabla u\vert^2 + s^4\la^4\va^4 \vert u\vert^2\\
\le& C\biggl\{ \frac{1}{\va} \left( \vert \ppp_tw\vert^2 
+ \sumij \vert \ppp_i\ppp_jw\vert^2 \right) 
+  \la^2\vert \nabla u\vert^2 + (\la^4+\va^2)\vert u\vert^2\biggr\}\\
+ &C(s^2\la^2\va\vert \nabla w\vert^2 + s^2\la^4\va^2\vert u\vert^2
+ s^4\la^4\va^3 \vert w\vert^2).
\end{align*}
In view of (2.1), we see
\begin{align*}
& \int_Q \left( \vert \ppp_tu\vert^2 + \sumij \vert \ppp_i\ppp_ju\vert^2
+ s^2\la^2\va^2 \vert \nabla u\vert^2 + s^4\la^4\va^4 \vert u\vert^2\right)
\weight dxdt \\
\le& C\int_Q \left\{ \frac{1}{\va} \left( \vert \ppp_tw\vert^2 
+ \sumij \vert \ppp_i\ppp_jw\vert^2 \right) 
+ s^2\la^2\va \vert \nabla w\vert^2 + s^4\la^4\va^3\vert w\vert^2\right\}
 \weight dxdt\\
+& C\int_Q (\la^2\vert \nabla u\vert^2 + (\la^4+\va^2)\vert u\vert^2 
+ s^2\la^4\va^2\vert u\vert^2) \weight dxdt\\
\le& C\int_Q s\va\vert F\vert^2 \weight dxdt
+ C\int_Q (s\va^3+s\la^4\va) \vert u\vert^2 \weight dxdt 
+ C\int_Q s\va\vert \nabla u\vert^2 \weight dxdt + C(s,\la)D(w)^2\\
+& C\int_Q (\la^2\vert \nabla u\vert^2 + (\la^4+\va^2)\vert u\vert^2
+ s^2\la^4\va^2 \vert u\vert^2)\weight dxdt.
\end{align*}
The second, the third and the fifth terms on the right-hand side can be 
absorbed into the left-hand side for large $s>0$.

Moreover, we can estimate the coefficient of $D(w)^2$ as follows:
\begin{align*}
& s\va^m \weight
= s\frac{e^{m\la\eta(x)}}{t^m(T-t)^m}
\exp\left( -2s \left( \frac{e^{2\la\EEEE}-e^{\la\eta(x)}}{t(T-t)}
\right) \right)\\
\le& C_1(\la)s\left( \frac{1}{t(T-t)}\right)^m
\exp\left( -2s\frac{C_2(\la)}{t(T-t)}\right) 
\le C_1(\la)s \max_{\xi \ge0} \xi^me^{-2sC_2(\la)\xi}\\
= & C_1(\la)s\left( \frac{m}{2sC_2(\la)}\right)^me^{-m} \quad 
\mbox{for each $m\in \N$},
\end{align*}
where the maximum is taken at $\xi = \frac{m}{2sC_2(\la)}$.
Thus the proof of Lemma 1 is complete.
$\blacksquare$
\\

Next
\\
{\bf Lemma 2 (Imanuvilov \cite{Ima}).}
\\
{\it
There exists a constant $\la_0>0$ such that for each $\la > \la_0$, we can find
a constant $s_0(\la) > 0$ satisfying: there exists a constant $C>0$ such that 
\begin{align*}
& \int_Q \biggl\{ \frac{1}{s\va}\left(
\vert \ppp_tv\vert^2 + \sumij \vert \ppp_i\ppp_jv\vert^2\right)
+ s\la^2\va\vert \nabla v\vert^2 + s^3\la^4\va^3\vert v\vert^2 \biggr\}
\weight dxdt\\
\le& C\int_Q \vert \ppp_tv - B(t)v\vert^2 \weight dxdt
+ C_1(\la,s)D(v)^2
\end{align*}
for all $s > s_0(\la)$ and $u\in H^{2,1}(Q)$ satisfying
$\ppp_{\nu_B}v = 0$ on $\ppp\OOO\times (0,T)$.
}
\\

Here the constant $C>0$ depends continuously on $M$: bound of the
coefficients and $\la$ but
independent of $s \ge s_0(\la)$, while $\la_0>0$ depends continuously
on $M$.
\\

Now, noting that we have the Carleman estimates both for 
$\ppp_t + A(t)$ and $\ppp_t - B(t)$ with the same weight $\weight$,
we derive the main Carleman estimate for the mean field game
equations (1.1).
Setting $F:= c_0v+F$, we apply Lemma 1 to obtain
$$
 \int_Q \left(\vert \ppp_tu\vert^2 + \sumij \vert \ppp_i\ppp_ju\vert^2
+ s^2\la^2\va^2\vert \nabla u\vert^2 + s^4\la^4\va^4\vert u\vert^2 \right)
\weight dxdt
$$
$$
\le C\int_Q s\va\vert v\vert^2 \weight dxdt 
+ C\int_Q s\va \vert F\vert^2 \weight dxdt
+ CD(u)^2.                                       \eqno{(2.3)}
$$
By Lemma 2, we obtain
$$
 \int_Q \biggl\{ \frac{1}{s\va}\left(
\vert \ppp_tv\vert^2 + \sumij \vert \ppp_i\ppp_jv\vert^2\right)
+ s\la^2\va\vert \nabla v\vert^2 + s^3\la^4\va^3\vert v\vert^2 \biggr\}
\weight dxdt
$$
$$
\le C\int_Q \sumij \vert \ppp_i\ppp_ju\vert^2 \weight dxdt
+ C\int_Q \vert G\vert^2 \weight dxdt + CD(v)^2.
                                                \eqno{(2.4)}
$$
Substituting (2.3) into the first term on the right-hand side of 
(2.4), we have
\begin{align*}
& \int_Q \biggl\{ \frac{1}{s\va}\left(
\vert \ppp_tv\vert^2 + \sumij \vert \ppp_i\ppp_jv\vert^2\right)
+ s\la^2\va\vert \nabla v\vert^2 + s^3\la^4\va^3\vert v\vert^2 \biggr\}
\weight dxdt\\
\le& C\int_Q s\va\vert v\vert^2 \weight dxdt
+ C\int_Q s\va\vert F\vert^2 \weight dxdt 
+ CD(u)^2 + C\int_Q \vert G\vert^2 \weight dxdt + CD(v)^2
\end{align*}
for all large $s, \la > 0$.
Hence, choosing $s>0$ and $\la > 0$ sufficiently large, we can absorb 
the first term on the right-hand side into the left-hand side, and we can 
obtain
\begin{align*}
& \int_Q \biggl\{ \frac{1}{s\va}\left(
\vert \ppp_tv\vert^2 + \sumij \vert \ppp_i\ppp_jv\vert^2\right)
+ s\la^2\va\vert \nabla v\vert^2 + s^3\la^4\va^3\vert v\vert^2 \biggr\}
\weight dxdt\\
\le& C\int_Q (s\va\vert F\vert^2 + \vert G\vert^2) \weight dxdt
+ C(D(u)^2 + D(v)^2)
\end{align*}
for all large $s, \la > 0$.  Adding with (2.3), we absorb the term
$\int_Q s\va \vert v\vert^2 \weight dxdt$ on the right-hand side into the 
left-hand side, so that we proved
\\
{\bf Theorem 3 (Carleman estimate for a generalized mean field game
equations)}.
\\
{\it
There exists a constant $\la_0>0$ such that for each $\la > \la_0$, 
we can find
a constant $s_0(\la) > 0$ satisfying: there exists a constant $C>0$ such that 
\begin{align*}
& \int_Q \biggl\{ \vert \ppp_tu\vert^2 + \sumij \vert \ppp_i\ppp_ju\vert^2
+ s^2\la^2\va^2\vert \nabla u\vert^2 + s^4\la^4\va^4\vert u\vert^2 \\
+ & \frac{1}{s\va}\left(
\vert \ppp_tv\vert^2 + \sumij \vert \ppp_i\ppp_jv\vert^2\right)
+ s\la^2\va\vert \nabla v\vert^2 + s^3\la^4\va^3\vert v\vert^2 \biggr\}
\weight dxdt\\
\le& C\int_Q (s\va\vert F\vert^2 + \vert G\vert^2) \weight dxdt
+ C(D(u)^2 + D(v)^2)
\end{align*}
for all $s > s_0(\la)$ and $u, v\in H^{2,1}(Q)$ satisfying
$\ppp_{\nu_A}u = \ppp_{\nu_B}v = 0$ on $\ppp\OOO\times (0,T)$.
}
\\

Here the constant $C>0$ depends continuously on $M$: bound of the
coefficients and $\la$ but
independent of $s \ge s_0(\la)$, while $\la_0>0$ depends continuously
on $M$.

Now we proceed to 
\\
{\bf Proof of Theorem 1.}
First we can reduce (1.4) to the case where $\ppp_{\nu}u = \ppp_{\nu}v = 0$ 
on $\ppp\OOO\times (0,T)$.
Indeed, by the extension theorem in Sobolev spaces (e.g., 
Theorem 9.4 (pp.41-42) in Lions and Magenes \cite{LM}), we can find
$\www{u}, \, \www{v} \in H^1(0,T;H^2(\OOO))$ such that 
$$
\www{u} = \www{v} = 0, \quad \ppp_{\nu_A}\www{u} = \rho, \quad
\ppp_{\nu_B}\www{v} = \sigma \quad \mbox{on $\ppp\OOO \times (0,T)$}
$$
and
$$
\Vert \www{u}\Vert_{H^1(0,T;H^2(\OOO))}
\le C\Vert \rho\Vert_{H^1(0,T;H^{\hhalf}(\ppp\OOO))}, \quad
\Vert \www{v}\Vert_{H^1(0,T;H^2(\OOO))}
\le C\Vert \sigma\Vert_{H^1(0,T;H^{\hhalf}(\ppp\OOO))}.  \eqno{(2.5)}
$$
Setting $U:= u - \www{u}$ and $V:= v - \www{v}$, we have
$$
\left\{ \begin{array}{rl}
& \ppp_tU + A(t)U = c_0V + \{F- (\ppp_t\www{u}+A_1(t)\www{u}-c_0\www{v})\},\\
& \ppp_tV - B(t)V = A_0(t)U 
    + \{ G - (\ppp_t\www{v} - B(t)\www{v} - A_0(t)\www{u})\} 
\quad \mbox{in $Q$},
\end{array}\right.
$$
and
$$
\ppp_{\nu_A}U = \ppp_{\nu_B}V = 0, \quad 
U=u, \quad V=v \quad \mbox{on $\ppp\OOO\times (0,T)$}.
$$
By (2.5), we see
\begin{align*}
& \Vert F- (\ppp_t\www{u}+A_1(t)\www{u}-c_0\www{v})\Vert_{L^2(Q)} 
\le \Vert F\Vert_{L^2(Q)} + \Vert \ppp_t\www{u}+A_1(t)\www{u}-c_0\www{v}
\Vert_{L^2(Q)}\\
\le&  \Vert F\Vert_{L^2(Q)} + C(\Vert \www{u}\Vert_{H^1(0,T;H^2(\OOO))}
+ \Vert \www{v}\Vert_{L^2(Q)})
\le \Vert F\Vert_{L^2(Q)} + C(\Vert \rho\Vert_{H^1(0,T;H^{\hhalf}(\ppp\OOO))}
+ \Vert \sigma\Vert_{H^1(0,T;H^{\hhalf}(\ppp\OOO))})
\end{align*}
and similarly 
$$
\Vert G- (\ppp_t\www{v} - A_2(t)\www{v}-A_0(t)\www{u})\Vert_{L^2(Q)}
\le \Vert G\Vert_{L^2(Q)} + C(\Vert \rho\Vert_{H^1(0,T;H^{\hhalf}(\ppp\OOO))}
+ \Vert \sigma\Vert_{H^1(0,T;H^{\hhalf}(\ppp\OOO))}).
$$
Therefore, it suffices to prove Theorem 1 by assuming that 
$\ppp_{\nu}u = \ppp_{\nu}v = 0$ on $\ppp\OOO \times (0,T)$.
\\

Now we can immediately complete the proof of Theorem 1.
Indeed, since 
$$
\alpha(x,t) \ge \frac{e^{\la\eta(x)} - e^{2\la\Vert \eta\Vert_{C(\ooo{\OOO})}} 
}{\ep (T-\ep)}
\ge \frac{-C_3}{\ep(T-\ep)} =: -C_4 < 0
$$
for all $x\in \ooo{\OOO}$ and $\ep \le t \le T-\ep$, we see that 
$$
e^{2s\alpha(x,t)} \ge e^{-2sC_4}, \quad x\in \ooo{\OOO}, \, 
\ep\le t \le T-\ep.
$$
Thus Theorem 2 completes the proof of Theorem 1.
$\blacksquare$
\section{Proof of Theorem 2}

We can execute the proof with a similar spirit to Theorem 3.1 
in Imanuvilov and Yamamoto \cite{IY98}, where we have to 
estimate extra second-order derivatives of $u$.  The key is Lemma 3 
stated below.

Without loss of generality, we can assume that $t_0 = \frac{T}{2}$ by 
scaling the time variable.
Indeed, we choose small $\delta > 0$ 
such that $0<t_0-\delta < t_0 < t_0 + \delta < T$.
Then we consider a change of the variables 
$t \mapsto \xi:= \frac{t-(t_0-\delta)}{2\delta}T$.  Then, the inverse problem 
over 
the time interval $(t_0-\delta,\, t_0+\delta)$ can be transformed to 
$(0,T)$ with $t_0= \frac{T}{2}$.  

We set $t_0 = \frac{T}{2}$.
\\
\vspace{0.2cm}
\\
{\bf Proof of Theorem 3.}
\\
{\it First Step: Carleman estimate for an integral operator
$\int^t_{t_0} \cdots d\xi$.}
\\
In this step, we prove 
\\
{\bf Lemma 3.}
\\
{\it
Let $p\ge 0$ and $w \in L^2(Q)$.  Then there exists a constant $C=C(p) > 0$ 
such that 
$$
\int_Q (s\va)^p \left\vert \int^t_{t_0} w(x,\xi) d\xi \right\vert^2 
\weight dxdt
\le C\int_Q (s\va)^{p-1}\la^{-1} \vert w(x,t)\vert^2 \weight dxdt
$$
for all large $s, \la > 0$.
}

A similar lemma is proved in Loreti, Sforza and Yamamoto 
\cite{LSY} where the weight is given by 
$e^{\la (\eta(x) - \beta (t-t_0)^2)}$.  For the case $p=0$, the works
Bukhgeim and Klibanov \cite{BK} and Klibanov \cite{Kli} firstly 
used this kind of weighted estimate
for the uniqueness of inverse problems.
\\
\\
{\bf Proof of Lemma 3.}
We divide the left-hand side and consider 
$$
\int^{T}_{t_0} (s\va)^p \left\vert \int^t_{t_0} w(x,\xi) d\xi \right\vert^2 
\weight dt
$$
and
$$
\int^{t_0}_0 (s\va)^p \left\vert \int^t_{t_0} w(x,\xi) d\xi \right\vert^2 
\weight dt.
$$
We set $h(x):= e^{2\la\Vert \eta\Vert_{C(\ooo{\OOO})}}
- e^{\la\eta(x)} > 0$.  Then we can write 
$\alpha(x,t) = \frac{-h(x)}{t(T-t)}$ in $Q$.
Moreover 
$$
\ppp_t(\weight) = 2s(\ppp_t\alpha)\weight
= 2s(-h(x))\frac{2(t-t_0)}{t^2(T-t)^2}\weight,
$$
and so 
$$
(t-t_0)(s\va)^p \weight = \frac{t^2(T-t)^2(s\va)^p}{4s(-h(x))}
\ppp_t(\weight).                           \eqno{(3.1)}
$$
We note
$$
\frac{t^2(T-t)^2(s\va)^p}{4s\vert h(x)\vert} \le \frac{C}{\la}(s\va)^{p-1},
\quad (x,t)\in Q.                            \eqno{(3.2)}
$$
\\
{\bf Proof of (3.2).}
\\
It suffices to prove
$$
\frac{t^2(T-t)^2}{\vert h(x)\vert} \le \frac{C}{\la\va}.
$$
Indeed, 
\begin{align*}
& \frac{t^2(T-t)^2\la\va}{\vert h(x)\vert}
= \frac{t(T-t)\la e^{\la\eta(x)}}
{e^{2\la\EEEE} - e^{\la\eta(x)}}
\le \frac{C\la e^{\la\eta(x)}}{e^{2\la\EEEE} - e^{\la\eta(x)}}\\
=& \frac{C\la}{e^{\la(2\EEEE-\eta(x))} - 1}
\le \frac{C\la}{e^{\la\EEEE} - 1}, \quad (x,t) \in Q.
\end{align*}
On the other hand, 
$$
\sup_{\la>0} \frac{\la}{e^{C_0\la} - 1} < \infty,
$$
where $C_0:= \EEEE$.  Therefore,
$$
\sup_{(x,t)\in Q} \frac{t^2(T-t)^2\la\va}{\vert h(x)\vert} < \infty,
$$
which proves (3.2).  $\blacksquare$
\\
Note 
$$
\ppp_t\va(x,t) = \frac{2(t-t_0)}{t^2(T-t)^2}e^{\la\eta(x)}.
                                                       \eqno{(3.3)}
$$
Let $t_0< t < T$.  Then by (3.1) and (3.2), the Cauchy-Schwarz inequality 
implies  
\begin{align*}
& \int^T_{t_0} (s\va)^p \left\vert \int^t_{t_0}w(x,\xi)d\xi 
\right\vert^2 \weight dt
\le \int^T_{t_0} (s\va)^p \left( \int^t_{t_0} \vert w(x,\xi)\vert^2 d\xi 
\right) (t-t_0) \weight dt\\
=& \int^T_{t_0} \left\{ \left( \int^t_{t_0} \vert w(x,\xi)\vert^2 d\xi\right)
\frac{(s\va)^pt^2(T-t)^2}{4s(-h(x))}\right\} \ppp_t(\weight) dt\\
=& \left[ \left( \int^t_{t_0} \vert w(x,\xi)\vert^2 d\xi\right)
\frac{(s\va)^pt^2(T-t)^2}{4s(-h(x))} e^{2s\alpha(x,t)}
\right]^{t=T}_{t=t_0}
+ \int^T_{t_0} \vert w(x,\xi)\vert^2
\frac{(s\va)^pt^2(T-t)^2}{4s\vert h(x)\vert}\weight dt\\
+& \int^T_{t_0} \left( \int^t_{t_0} \vert w(x,\xi)\vert^2 d\xi\right)
\frac{(s\va)^p t(t-2t_0)(t-t_0)}{s\vert h(x)\vert} \weight dt\\
+ & \int^T_{t_0} \left\{ \int^t_{t_0} \vert w(x,\xi)\vert^2 d\xi\right) 
\frac{p}{2}(s\va)^{p-1} \frac{(t-t_0)e^{\la\eta(x)}}{\vert h(x)\vert}
\weight dt\\
\le& C\int^T_{t_0} \vert w(x,t)\vert^2 \frac{1}{\la}(s\va)^{p-1}
\weight dt
+ C\int^T_{t_0}\left(\int^t_{t_0} \vert w(x,\xi)\vert^2 d\xi\right)
(s\va)^{p-1}\weight dt.
\end{align*}
Here we used also (3.3), $\alpha(x,T) = -\infty$,
$$
\frac{t(t-2t_0)(t-t_0)}{\vert h(x)\vert} \le 0
$$
and
\begin{align*}
& \frac{e^{\la\eta(x)}}{\vert h(x)\vert}
\le C\frac{e^{\la\eta(x)}}{e^{2\la\EEEE} - e^{\la\eta(x)}}
= \frac{C}{e^{\la(2\EEEE - \eta(x))} - 1}\\
\le & \frac{C}{e^{\la(2\EEEE - \EEEE)} - 1} = \frac{C}{e^{\la\EEEE} - 1}
\quad \mbox{for $\la > 1$},
\end{align*}
which implies 
$$
\sup_{\la>1}\sup_{(x,t)\in Q} \frac{e^{\la\eta(x)}}{\vert h(x)\vert}
< \infty.                \eqno{(3.4)}
$$
Hence,
\begin{align*}
& \int_{\OOO} \int^T_{t_0} (s\va)^p \left\vert 
\int^t_{t_0} \vert w(x,\xi)\vert^2 d\xi \right\vert \weight dtdx\\
\le& C\int_{\OOO} \int^T_{t_0} s^{p-1}\la^{-1}\va^{p-1} 
\vert w(x,t)\vert^2 \weight dtdx
+ C\int_{\OOO} \int^T_{t_0} (s\va)^{p-1} 
\left\vert \int^t_{t_0} \vert w(x,\xi)\vert^2 d\xi \right\vert^2 
\weight dtdx
\end{align*}
for all large $s, \la >0$.
Choosing $s>0$ and $\la > 1$ sufficiently large, we can absorb the 
second term on the right-hand side into the left-hand side.
Thus
$$
\int_{\OOO} \int^T_{t_0} (s\va)^p \left\vert 
\int^t_{t_0} \vert w(x,\xi)\vert^2 d\xi \right\vert \weight dtdx
\le C\int_{\OOO} \int^T_{t_0} s^{p-1}\la^{-1}\va^{p-1} 
\vert w(x,t)\vert^2 \weight dtdx.                \eqno{(3.5)}
$$

Similarly 
\begin{align*}
& \int^{t_0}_0 (s\va)^p \left\vert \int^t_{t_0}w(x,\xi)d\xi 
\right\vert^2 \weight dt
= \int^{t_0}_0 (s\va)^p \left\vert \int^{t_0}_tw(x,\xi)d\xi 
\right\vert^2 \weight dt\\
\le & \int^{t_0}_0 (s\va)^p \left( \int^{t_0}_t \vert w(x,\xi)\vert^2 d\xi 
\right) (t_0-t) \weight dt\\
=& \int^{t_0}_0 \left\{ \left( \int^{t_0}_t \vert w(x,\xi)\vert^2 d\xi\right)
\frac{(s\va)^pt^2(T-t)^2}{4s\vert h(x)\vert}\right\} \ppp_t(\weight) dt\\
=& \left[ \left( \int^{t_0}_t \vert w(x,\xi)\vert^2 d\xi\right)
\frac{(s\va)^pt^2(T-t)^2}{4s\vert h(x)\vert} e^{2s\alpha(x,t)}
\right]^{t=t_0}_{t=0}
+ \int^{t_0}_0 \vert w(x,t)\vert^2
\frac{(s\va)^pt^2(T-t)^2}{4s\vert h(x)\vert}\weight dt\\
-& \int^{t_0}_0 \left( \int^{t_0}_t  \vert w(x,\xi)\vert^2 d\xi\right)
\frac{(s\va)^p t(t-2t_0)(t-t_0)}{s\vert h(x)\vert} \weight dt\\
- & \int^{t_0}_0 \left\{ \int^{t_0}_t \vert w(x,\xi)\vert^2 d\xi\right) 
\frac{p}{2}(s\va)^{p-1} \frac{(t-t_0)e^{\la\eta(x)}}{\vert h(x)\vert}
\weight dt\\
\le& C\int^{t_0}_0 \vert w(x,t)\vert^2 \frac{t^2(T-t)^2}{4s\vert h(x)\vert}
(s\va)^p \weight dt
+ C\int^{t_0}_0 \left( \int^{t_0}_t \vert w(x,\xi)\vert^2 d\xi\right)
(s\va)^{p-1}\weight dt.
\end{align*}
Here we used (3.3), (3.4), $\alpha(x,0) = -\infty$ and
$$
\frac{t(t-2t_0)(t-t_0)}{s\vert h(x)\vert} \ge 0 \quad 
\mbox{for all $x\in \OOO$ and $0<t<t_0$}.
$$
Therefore (3.2) verifies
\begin{align*}
& \int^{t_0}_0 (s\va)^p \left\vert \int^t_{t_0}
\vert w(x,\xi)\vert^2 d\xi\right\vert \weight dt\\
\le& C\int^{t_0}_0 s^{p-1}\la^{-1}\va^{p-1} \vert w(x,t)\vert^2 
\weight dt
+ C\int^{t_0}_0 (s\va)^{p-1} \left( \int^{t_0}_t \vert w(x,\xi)\vert^2 d\xi
\right) dt.
\end{align*}
Choosing $s, \la > 0$ large, we can reach 
$$
\int^{t_0}_0 (s\va)^p \left\vert \int^t_{t_0}
\vert w(x,\xi)\vert^2 d\xi\right\vert \weight dt
\le C\int^{t_0}_0 s^{p-1}\la^{-1}\va^{p-1} \vert w(x,t)\vert^2 
\weight dt.
$$
Thus by combining with (3.5), the proof of Lemma 3 is complete.
$\blacksquare$
\\
\vspace{0.2cm}
\\
{\it Second Step: Carleman estimate for $(\ppp_tu, \ppp_tv)$.}
\\
Setting $y:= \ppp_tu$ and $z:= \ppp_tz$, we have
$$
\left\{ \begin{array}{rl}
& \ppp_ty + A(t)y = c_0z + (\ppp_tc_0)v - (\ppp_tA(t))u + \ppp_tF,\\
& \ppp_tz - B(t)z = A_0(t)y + (\ppp_tA_0(t))u - (\ppp_tB(t))v + \ppp_tG
\end{array}\right.
$$
in $Q$ and
$$
\ppp_{\nu_A}y = \ppp_{\nu_B}z = 0 \quad \mbox{on $\ppp\OOO \times (0,T)$}.
$$

Application of Theorem 3 yields
\begin{align*}
& \int_Q \biggl\{ \vert \ppp_t^2u\vert^2 
+ \sumij \vert \ppp_i\ppp_j\ppp_tu\vert^2
+ s^2\la^2\va^2\vert \nabla \ppp_tu\vert^2 
+ s^4\la^4\va^4\vert \ppp_tu\vert^2 \\
+ & \frac{1}{s\va}\left(
\vert \ppp_t^2v\vert^2 + \sumij \vert \ppp_i\ppp_j\ppp_tv\vert^2\right)
+ s\la^2\va\vert \nabla \ppp_tv\vert^2 + s^3\la^4\va^3\vert \ppp_tv\vert^2
\biggr\} \weight dxdt\\
\le& C\int_Q s\va\left( \vert v\vert^2 + \sumij \vert \ppp_i\ppp_ju\vert^2
\right) \weight dxdt
+ C\int_Q \left( \sumij \vert \ppp_i\ppp_ju\vert^2 
+ \sumij \vert \ppp_i\ppp_jv\vert^2 \right) \weight dxdt
\end{align*}
$$
+ C\int_Q (s\va\vert \ppp_tF\vert^2 + \vert \ppp_tG\vert^2) \weight dxdt 
+ C(D(\ppp_tu)^2 + D(\ppp_tv))^2)                \eqno{(3.6)}
$$
for all large $s, \la > 0$.

It is necessary to estimate the terms of the second-order derivatives of
$u$ and $v$ on the right-hand side of (3.6).  For it, we will apply Lemma 3.
First, since
$$
\ppp_i\ppp_ju(x,t) = \int^t_{t_0} \ppp_i\ppp_j\ppp_tu(x,\xi) d\xi
+ \ppp_i\ppp_ju_0(x), \quad (x,t) \in Q,
$$
we estimate
\begin{align*}
& \int_Q s\va\vert \ppp_i\ppp_ju(x,t)\vert^2 \weight dxdt
\\
\le& C\int_Q s\va \left\vert \int^t_{t_0} \ppp_i\ppp_j\ppp_tu(x,\xi) d\xi
\right\vert^2 \weight dxdt 
+ C\int_Q s\va \vert \ppp_i\ppp_ju_0(x)\vert^2 \weight dxdt.
\end{align*}
Applying Lemma 3 with $p=1$, we obtain
$$
 \int_Q s\va\vert \ppp_i\ppp_ju(x,t)\vert^2 \weight dxdt
\le \frac{C}{\la} \int_Q \vert \ppp_i\ppp_j\ppp_tu(x,t)\vert^2 \weight dxdt
+ C\Vert u_0\Vert^2_{H^2(\OOO)}.
$$
Similarly we can verify
$$
 \int_Q \vert \ppp_i\ppp_jv(x,t)\vert^2 \weight dxdt
\le \frac{C}{\la} \int_Q \frac{1}{s\va}\vert \ppp_i\ppp_j\ppp_tv(x,t)\vert^2 
\weight dxdt + C\Vert v_0\Vert^2_{H^2(\OOO)}.
$$
Therefore, choosing $\la>0$ sufficiently large, we can absorb all the terms
$\vert \ppp_i\ppp_ju\vert^2$ and $\vert \ppp_i\ppp_jv\vert^2$ on the 
right-hand side into the left-hand side of (3.6).
By Theorem 3, we can absorb $\int_Q s\va \vert v\vert^2 \weight dxdt$ and 
we can reach
\\
{\bf Lemma 4 (Carleman estimate for $(\ppp_tu, \, \ppp_tv)$).}
\\
{\it
\begin{align*}
& \int_Q \biggl\{ \vert \ppp_t^2u\vert^2 
+ \sumij \vert \ppp_i\ppp_j\ppp_tu\vert^2
+ s^2\la^2\va^2\vert \nabla \ppp_tu\vert^2 
+ s^4\la^4\va^4\vert \ppp_tu\vert^2 \\
+ & \frac{1}{s\va}\left(
\vert \ppp_t^2v\vert^2 + \sumij \vert \ppp_i\ppp_j\ppp_tv\vert^2\right)
+ s\la^2\va\vert \nabla \ppp_tv\vert^2 
+ s^3\la^4\va^3\vert \ppp_tv\vert^2 \biggr\}
\weight dxdt\\
\le& C\int_Q (s\va(\vert \ppp_tF\vert^2 + \vert F\vert^2) 
+ \vert \ppp_tG\vert^2 + \vert G\vert^2) \weight dxdt
+ CD_0^2
\end{align*}
for all large $s, \la > 0$.
}
\\
Here and henceforth, we set 
$$
D_0^2:= D(u)^2 + D(v)^2 + D(\ppp_tu)^2 + D(\ppp_tv)^2
+ \Vert u_0\Vert^2_{H^2(\OOO)} + \Vert v_0\Vert^2_{H^2(\OOO)}.
$$

We apply Lemma 4 to the system in $\ppp_tu, \ppp_tv$ with 
$F(x,t) = q_1(x,t)f(x)$ and $G(x,t) = q_2(x,t)g(x)$, and, fixing 
$\la > 0$ large, we reach
$$
 \int_Q \biggl\{ \vert \ppp_t^2u\vert^2 
+ s^4\va^4\vert \ppp_tu\vert^2 
+ \frac{1}{s\va}\vert \ppp_t^2v\vert^2 + s^3\va^3\vert \ppp_tv\vert^2 \biggr\}
\weight dxdt
$$
$$
\le C\int_Q (s\va\vert f\vert^2 + \vert g\vert^2) \weight dxdt
+ CD_0^2
                                              \eqno{(3.7)}
$$
for all large $s > 0$.
\\
{\it Third Step: Completion of the proof.}
\\
We note that $\vert \ppp_t\alpha \vert \le C\vert \va\vert^2$ in $Q$
with fixed large $\la>0$.
We have
\begin{align*}
& \int_{\OOO} s\va(x,t_0) \vert \ppp_tu(x,t_0)\vert^2
e^{2s\alpha(x,t_0)} dx 
= \int^{t_0}_0 \ppp_t\left( \int_{\OOO} s\va\vert \ppp_tu(x,t)\vert^2
e^{2s\alpha(x,t)} dx \right) dt \\
=& 2\int^{t_0}_0 \int_{\OOO} (2s\va (\ppp_tu)(\ppp_t^2u)
+ s(\ppp_t\va)\vert \ppp_tu\vert^2 
+ 2s\va \times s(\ppp_t\alpha)\vert \ppp_tu\vert^2) \weight dxdt\\
\le& C\int_Q (s\va\vert \ppp_tu\vert \vert \ppp_t^2u\vert
+ s\va^2 \vert \ppp_tu\vert^2 + s^2\va^3\vert \ppp_tu\vert^2) \weight dxdt \\
= & C\int_Q ((s^{-\hhalf}\vert \ppp_t^2u\vert)
(s^{\frac{3}{2}}\va \vert \ppp_tu\vert) 
+ (s\va^2 + s^2\va^3)\vert \ppp_tu\vert^2) \weight dxdt \\
\le & C\int_Q (s^{-1}\vert \ppp_t^2u\vert^2
+ (s^3\va^2 + s\va^2 + s^2\va^3)\vert \ppp_tu\vert^2) \weight dxdt
\le \frac{C}{s}\int_Q (\vert \ppp_t^2u\vert^2 + s^4\va^4\vert \ppp_tu\vert^2)
 \weight dxdt.
\end{align*}
Hence, by (3.7), we have
$$
 \int_{\OOO} s\va(x,t_0) \vert \ppp_tu(x,t_0)\vert^2
e^{2s\alpha(x,t_0)} dx 
\le C\int_Q \left(\va \vert f\vert^2 + \frac{1}{s}\vert g\vert^2\right) 
\weight dxdt + CD_0^2                     \eqno{(3.8)}
$$
for all large $s>0$.

Next we can similarly estimate
\begin{align*}
& \int_{\OOO} \vert \ppp_tv(x,t_0)\vert^2
e^{2s\alpha(x,t_0)} dx 
= \int^{t_0}_0 \ppp_t\left( \int_{\OOO} \vert \ppp_tv(x,t)\vert^2
e^{2s\alpha(x,t)} dx \right) dt \\
=& 2\int^{t_0}_0 \int_{\OOO} (2 (\ppp_tv)(\ppp_t^2v)
+ 2s(\ppp_t\alpha)\vert \ppp_tv\vert^2) \weight dxdt\\
\le& C\int_Q (\vert \ppp_tv\vert \vert \ppp_t^2v\vert
+ s\va^2 \vert \ppp_tv\vert^2) \weight dxdt 
= C\int_Q \left\{
\left( \frac{1}{s^{\frac{3}{4}}\va^{\hhalf}} \vert \ppp_t^2v\vert\right)
\left( s^{\frac{3}{4}}\va^{\hhalf} \vert \ppp_tv\vert\right) 
+ s\va^2 \vert \ppp_tv\vert^2) \weight \right\}dxdt \\
\le & \frac{C}{\sqrt{s}}\int_Q \left( \frac{1}{s\va} \vert \ppp_t^2v\vert^2
+ s^2\va\vert \ppp_tv\vert^2 + s^{\frac{3}{2}} \va^2\vert \ppp_tv\vert^2
\right) \weight dxdt\\
\le & \frac{C}{\sqrt{s}}\int_Q \left(
\frac{1}{s\va}\vert \ppp_t^2v\vert^2 + s^3\va^3\vert \ppp_tv\vert^2\right)
 \weight dxdt.
\end{align*}
The estimate (3.7) implies
$$
 \int_{\OOO} \vert \ppp_tv(x,t_0)\vert^2
e^{2s\alpha(x,t_0)} dx 
\le C\int_Q \left(\sqrt{s}\va \vert f\vert^2 
+ \frac{1}{\sqrt{s}}\vert g\vert^2\right) \weight dxdt 
+ CD_0^2                     \eqno{(3.9)}
$$
for all large $s>0$.

On the other hand, we have
$$
\left\{ \begin{array}{rl}
& f(x) = \frac{1}{q_1(x,t_0)}\ppp_tu(x,t_0)
+ \frac{1}{q_1(x,t_0)}(A_1(t_0)u_0 - c_0(x,t_0)v_0),\\
& g(x) = \frac{1}{q_2(x,t_0)}\ppp_tv(x,t_0)
+ \frac{1}{q_2(x,t_0)}(-A_2(t_0)v_0 - A_0(t_0)u_0), \quad x\in \OOO
\end{array}\right.
$$
and
$$
\left\{ \begin{array}{rl}
& s\va(x_0)\vert f(x)\vert^2
\le Cs\va(x,t_0)\vert \ppp_tu(x,t_0)\vert^2
+ Cs\va(x,t_0) \left( \sum_{\vert \gamma\vert\le 2}
\vert \ppp_x^{\gamma}u_0(x)\vert^2 + \vert v_0(x)\vert^2\right),\\
& \vert g(x)\vert^2
\le C\vert \ppp_tu(x,t_0)\vert^2
+ C\sum_{\vert \gamma\vert\le 2}(\vert \ppp_x^{\gamma}u_0(x)\vert^2 
+ \vert \ppp_x^{\gamma} v_0(x)\vert^2).
\end{array}\right.
$$
Hence (3.8) and (3.9) yield
\begin{align*}
& \int_{\OOO} (s\va(x_0,t)\vert f(x)\vert^2 + \vert g(x)\vert^2)
e^{2s\alpha(x,t_0)} dx\\
\le& C\int_{\OOO} (s\va(x,t_0)\vert \ppp_tu(x,t_0)\vert^2
+ \vert \ppp_tv(x,t_0)\vert^2) \weight dxdt
+ C(\Vert u_0\Vert^2_{H^2(\OOO)} + \Vert v_0\Vert^2_{H^2(\OOO)})\\
\le& C\int_Q \left( \va\vert f\vert^2 + \frac{1}{s}\vert g\vert^2\right)
\weight dxdt 
+ C\int_Q \left( \sqrt{s}\va\vert f\vert^2 
+ \frac{1}{\sqrt{s}}\vert g\vert^2\right)\weight dxdt + CD_0^2,
\end{align*}
that is,
$$
 \int_{\OOO} (s\va(x,t_0)\vert f(x)\vert^2 + \vert g(x)\vert^2)
e^{2s\alpha(x,t_0)} dx
$$
$$
\le C\int_Q \left( \sqrt{s}\va(x,t)\vert f(x)\vert^2 
+ \frac{1}{\sqrt{s}}\vert g(x)\vert^2\right)e^{2s\alpha(x,t)} dxdt 
+ CD_0^2                      \eqno{(3.10)}
$$
for all large $s>0$.
Here, setting $\ell(t):= t(T-t)$, we see
\begin{align*}
& \va(x,t)\vert f(x)\vert^2 e^{2s\alpha(x,t)}
= \va(x,t_0)\vert f(x)\vert^2 e^{2s\alpha(x,t_0)}
\times \frac{\va(x,t)}{\va(x,t_0)}e^{2s(\alpha(x,t)-\alpha(x,t_0))}\\
=& \va(x,t_0)\vert f(x)\vert^2 e^{2s\alpha(x,t_0)}
\times \frac{t_0^2}{4e^{\la\eta(x)}}\frac{e^{\la\eta(x)}}{\ell(t)}
e^{-2sh(x)\left( \frac{1}{\ell(t)} - \frac{1}{\ell(t_0)}\right)}.
\end{align*}
Here we set 
$$
h(x):= e^{2\la\Vert \eta\Vert_{C(\ooo{\OOO})}}
- e^{\la\eta(x)}> 0, \quad
C_1:= e^{2\la\Vert \eta\Vert_{C(\ooo{\OOO})}}
- e^{\la\Vert \eta\Vert_{C(\ooo{\OOO})}} > 0.
$$
Therefore,
\begin{align*}
& \int_Q \va(x,t)\vert f(x)\vert^2 e^{2s\alpha(x,t)} dxdt
\le C\int_Q \va(x,t_0)\vert f(x)\vert^2 e^{2s\alpha(x,t_0)}
\frac{1}{\ell(t)}
e^{-2sh(x)\left( \frac{1}{\ell(t)} - \frac{1}{\ell(t_0)}\right)} dxdt\\
\le& C\int_{\OOO} \va(x,t_0)\vert f(x)\vert^2 e^{2s\alpha(x,t_0)}
\left( \int^T_0 \frac{1}{\ell(t)}
e^{-2sC_1\left( \frac{1}{\ell(t)} - \frac{1}{\ell(t_0)}\right)}dt
\right) dx.
\end{align*}
We will estimate 
$\int^T_0 \frac{1}{\ell(t)}
e^{-2sC_1\left( \frac{1}{\ell(t)} - \frac{1}{\ell(t_0)}\right)}dt$.
Indeed 
$$
\lim_{s\to \infty} 
\frac{1}{\ell(t)}
e^{-2sC_1\left( \frac{1}{\ell(t)} - \frac{1}{\ell(t_0)}\right)} = 0
$$
for each fixed $t \in (0,T) \setminus \{t_0\}$.  Next, since 
$$
\frac{1}{\ell(t)}e^{-2sC_1\left( \frac{1}{\ell(t)} - \frac{1}{\ell(t_0)}\right)}\le \frac{1}{\ell(t)}
e^{-2C_1\left( \frac{1}{\ell(t)} - \frac{1}{\ell(t_0)}\right)}
$$
for $s \ge 1$ and by $\ell(t) \le \ell(t_0)$, we have
\begin{align*}
& \sup_{s\ge 1}\sup_{0< t< T} \frac{1}{\ell(t)}
\exp\left( -2sC_1\left( \frac{1}{\ell(t)} - \frac{1}{\ell(t_0)} \right)\right)
\le \sup_{0<t<T} \left( \frac{1}{\ell(t)}
e^{-\frac{2C_1}{\ell(t)}}\right) e^{\frac{2C_1}{\ell(t_0)}}\\ 
\le &\sup_{\xi>0} (\xi e^{-2C_1\xi}) e^{\frac{2C_1}{\ell(t_0)}} < \infty.
\end{align*}
Consequently, the Lebesgue convergence theorem yields
$$
\int^T_0 \frac{1}{\ell(t)}
e^{-2sC_1\left( \frac{1}{\ell(t)} - \frac{1}{\ell(t_0)}\right)} dt
= o(1) \quad \mbox{as $s \to \infty$}.
$$
Hence,
$$
\int_Q \va(x,t)\vert f(x)\vert^2 e^{2s\alpha(x,t)} dxdt
= o(1)\int_{\OOO} \va(x,t_0)\vert f(x)\vert^2 e^{2s\alpha(x,t_0)} dx
$$
as $s \to \infty$.  By (3.10), using $\alpha(x,t) \le \alpha(x,t_0)$ for
$(x,t) \in Q$, we obtain
$$
\int_{\OOO} (s\va(x,t_0)\vert f(x)\vert^2 + \vert g(x)\vert^2)
e^{2s\alpha(x,t_0)} dx
\le C\int_{\OOO} (o(\sqrt{s})\va(x,t_0)\vert f(x)\vert^2 
+ \frac{1}{\sqrt{s}}\vert g(x)\vert^2)e^{2s\alpha(x,t_0)} dx
+ CD_0^2
$$
for all large $s>0$.
Choosing $s>0$ sufficiently large, we absorb the first term on the right-hand 
side into the left-hand side, we reach 
$$
\int_{\OOO} (s\va(x_0)\vert f(x)\vert^2 + \vert g(x)\vert^2)
e^{2s\alpha(x,t_0)} dx
\le CD_0^2,
$$
which completes the proof of Theorem 2.
$\blacksquare$
\section{State determination for the nonlinear mean field game system}

For the not-linearlized mean field game system (1.1), we discuss the state
determination problem.  In this article, we do not consider inverse source 
problems and inverse coefficient problems of determining spatially varying
factors of the coefficients such as $\kappa(x,t)$, 
and we postpone them to a future work.  

In (1.1), we assume that 
$$
a \in C^1(\ooo{Q}) \cap L^{\infty}(0,T;W^{2,\infty}(\OOO)),\, > 0 
\quad \mbox{on $\ooo{Q}$}, \quad
p\in L^{\infty}(Q), \quad \kappa \in L^{\infty}(0,T;W^{1,\infty}(\OOO)).
                                               \eqno{(4.1)}
$$
Then we can prove
\\
{\bf Theorem 4.}
\\
{\it
For $k=1,2$, let $(u_k,v_k) \in H^{2,1}(Q)$ satisfy 
$$
\left\{ \begin{array}{rl}
& \ppp_tu_k + a(x,t)\Delta u_k(x,t) 
- \hhalf\kappa(x,t)\vert \nabla u_k(x,t)\vert^2 + p(x,t)v_k = F_k(x,t), \\
& \ppp_tv_k - \Delta (a(x,t) v_k(x,t)) 
- \ddd(\kappa(x,t)v_k(x,t)\nabla u_k(x,t)) = G_k(x,t), \quad (x,t) \in Q
\end{array}\right.
                                                          \eqno{(4.2)}
$$
and
$$
\nabla u_k\cdot \nu = \nabla v_k\cdot \nu = 0 
\quad \mbox{on $\ppp\OOO\times (0,T)$}.
$$
We further assume
$$
\Vert u_k\Vert_{L^{\infty}(0,T;W^{2,\infty}(\OOO))}
+ \Vert v_k\Vert_{L^{\infty}(0,T;W^{1,\infty}(\OOO))} \le M_1,  \eqno{(4.3)}
$$
where $M_1>0$ is an arbitrarily chosen constant.
Then, for arbitrarily given $\ep > 0$, we can find a constant $C_{\ep} >0$ 
such that 
\begin{align*}
& \Vert u_1-u_2\Vert_{\HHHH(\OOO\times (\ep,T-\ep))}
+ \Vert v_1-v_2\Vert_{\HHHH(\OOO\times (\ep,T-\ep))}
\le C_{\ep}(\Vert F_1-F_2\Vert_{L^2(Q)} + \Vert G_1-G_2\Vert_{L^2(Q)}\\
+ & \Vert u_1-u_2\Vert_{H^1(0,T;L^2(\Gamma))}
+ \Vert \nabla (u_1-u_2)\Vert_{L^2(\Gamma \times (0,T))}
+ \Vert v_1-v_2\Vert_{H^1(0,T;L^2(\Gamma))}
+ \Vert \nabla (v_1-v_2)\Vert_{L^2(\Gamma \times (0,T))}).
\end{align*}
}

Here the constant $C_{\ep}>0$ depends on $\ep, M_1 > 0$ and 
the coefficients $a, p, \kappa$.
\\
{\bf Proof of Theorem 4.}
\\
We take the difference $u:= u_1-u_2$ and $v:= v_1 - v_2$.  Then 
$$
\left\{ \begin{array}{rl}
& \ppp_tu + a(x,t)\Delta u(x,t) 
- \hhalf\kappa(x,t)((\nabla u_1 + \nabla u_2)\cdot\nabla u) + p(x,t)v 
= F_1-F_2, \\
& \ppp_tv - a(x,t)\Delta v(x,t) -2\nabla a(x,t)\cdot \nabla v(x,t)
- v(x,t)\Delta a(x,t)\\ 
- & (\kappa\nabla u_1\cdot \nabla v 
+ (\nabla\kappa \cdot \nabla u_1 + \kappa\Delta u_1)v 
+ v_2\kappa \Delta u + \nabla(\kappa v_2)\cdot \nabla u)
= G_1 - G_2, \quad (x,t) \in Q
\end{array}\right.
$$
By (4.1) and (4.3), we can find a constant $M_2 > 0$ depending on 
$M_1$ such that 
\begin{align*}
& \left\Vert \sum_{k=1}^2 \kappa\nabla u_k\right\Vert_{L^{\infty}(Q)}
+ \Vert\nabla \kappa \cdot \nabla u_1 + \kappa\Delta u_1\Vert_{L^{\infty}(Q)}\\
+& \Vert \kappa v_2\Vert_{L^{\infty}(Q)}  
+ \Vert\nabla (\kappa v_2)\Vert_{L^{\infty}(Q)} \le M_2.
\end{align*}
Therefore, we can satisfy the conditions of Theorem 3 where the constants
$C>0$, $\la$ and $s_0(\la)$ depend only on $M_2, a(x,t), \kappa(x,t), p(x,t)$, 
but independent of choices $u_1,v_1,u_2, v_2$.
Thus the proof of Theorem 4 follows directly from Theorem 3.
$\blacksquare$
\\
\vspace{0.2cm}
\\
{\bf Acknowledgments.}
The work was supported by Grant-in-Aid for Scientific Research (A) 20H00117 
of Japan Society for the Promotion of Science.


\begin{thebibliography}{99} %

\bibitem{ACDPS}
Y. Achdou, P. Cardaliaguet, F. Delarue, A. Porretta and F. Santambrogio, 
{\it Mean field games},
Cetraro, Italy 2019, Lecture Notes in Mathematics, C.I.M.E. 
Foundation Subseries, Volume 2281, Springer, 2019.

\bibitem{BK}
A.L. Bukhgeim and M.V. Klibanov, 
Global uniqueness of a class of multidimensional inverse problems,
Sov. Math.-Dokl. {\bf 24} (1981) 244-247.

\bibitem{CCP}
P. Cardaliaguet, M. Cirant and A. Porretta, 
Splitting methods and short time existence for the master equations in mean 
field games, J. Eur. Math. Soc. (2022),
DOI 10.4171/JEMS/1227

\bibitem{FI}
A. V. Fursikov and O. Y. Imanuvilov,  {\it Controllability of Evolution 
Equations}, Lecture Notes Series vol 34, 1996, Seoul National University.

\bibitem{Ima}
O.Y. Imanuvilov, Controllability of parabolic equations, Sbornik Math. 
{\bf 186} (1995) 879-900.

\bibitem{IY98}
O.Y. Imanuvilov and M. Yamamoto,
Lipschitz stability in inverse parabolic problems by the Carleman
estimate, Inverse Problems {\bf 14} (1998) 1229-1245.

\bibitem{Kli}
M.V. Klibanov,  Inverse problems and Carleman estimates, Inverse Problems 
{\bf 8} (1992) 575-596.

\bibitem{Kl23}
M. V. Klibanov, 
The mean field games system: Carleman estimates, Lipschitz stability and 
uniqueness, preprint arXiv:2303.03928

\bibitem{KlAv}
M. V. Klibanov and Y. Averboukh,
Lipschitz stability estimate and uniqueness in the retrospective analysis 
for the mean field games system via two Carleman estimates,
preprint arXiv:2302.10709


\bibitem{KLL1}
M. V. Klibanov, J. Li and H. Liu,
H\"older stability and uniqueness for the mean field games system 
via Carleman estimates, preprint, arXiv:2304.00646

\bibitem{KLL2}
M. V. Klibanov, J. Li and H. Liu,
On the mean field games system with the lateral Cauchy data via 
Carleman estimates, preprint, arXiv:2303.07556 

\bibitem{LL}
J.-M. Lasry and P.-L. Lions, Mean field games, 
Japanese Journal of Mathematics, {\bf 2} (2007) 229-260.

\bibitem{LM}
J.-L. Lions and E. Magenes, {\it
Non-homogeneous Boundary Value Problems and Applications},
vol.1, Springer-Verlag, Berlin, 1972.

\bibitem{LY}
H. Liu and M. Yamamoto, Stability in determination of states for the 
mean field game equations, preprint.

\bibitem{LZ1}
H. Liu and S. Zhang, On an inverse boundary problem for mean field games,
preprint, arXiv:2212.09110 

\bibitem{LZ2}
H. Liu and S. Zhang, Simultaneously recovering running cost and Hamiltonian 
in mean field games system, preprint,
arXiv:2303.13096

\bibitem{LSY}
P. Loreti, D. Sforza and M. Yamamoto, Carleman estimate and 
application to an inverse source problem for a viscoelasticity model 
in anisotropic case, Inverse Problems {\bf 33} (2017) 125014, 28 pp. 

\bibitem{Y09}
M. Yamamoto, Carleman estimates for parabolic equations and applications, 
Inverse Problems {\bf 25} (2009) 123013


\end{thebibliography}
\end{document}